\documentclass{amsproc}%
\usepackage{amsmath}
\usepackage{amsfonts}
\usepackage{amssymb}
\usepackage{graphicx}%
\theoremstyle{plain}

\newtheorem{corollary}{Corollary}

\newtheorem{definition}{Definition}

\newtheorem{lemma}{Lemma}

\newtheorem{theorem}{Theorem}
\numberwithin{equation}{section}

\begin{document}
\title[Removing Coincidences of Maps]{Removing Coincidences of Maps Between Manifolds of Different Dimensions.}
\author{Peter Saveliev}
\address{Marshall University, Huntington, WV 25755, USA}
\email{saveliev@member.ams.org}
\urladdr{http://saveliev.net}
\thanks{}

\begin{abstract}
We consider sufficient conditions of local removability of coincidences of
maps $f,g:N\rightarrow M,$ where $M,N$ are manifolds with dimensions $\dim
N\geq\dim M.$ The coincidence index is the only obstruction to the
removability for maps with fibers either acyclic or homeomorphic to spheres of
certain dimensions$.$ We also address the normalization property of the index
and coincidence-producing maps.

\end{abstract}
\date{}
\maketitle

\section{Introduction}

Let\ $N^{n+m}$ and $M^{n}$ be orientable compact smooth manifolds (possibly
with boundaries $\partial N,\partial M)$, $n\geq2,$ and suppose
$f,g:N\rightarrow M$ are maps. We shall call $m$ the \textit{codimension}. The
\textit{coincidence set }is a compact subset of $N$ defined by\textit{
}$Coin(f,g)=\{x\in N:f(x)=g(x)\}.$

The Coincidence Problem asks what can be said about the coincidence set. When
$m=0,$ the main tools for studying the problem is the \textit{Lefschetz
number} $L(f,g)$ defined as the alternating sum of traces of certain
endomorphism on the (co)homology group of $M$. The famous Lefschetz
coincidence theorem provides a sufficient condition for the existence of
coincidences: $L(f,g)\neq0\Longrightarrow Coin(f,g)\neq\emptyset.$ Under some
circumstances the converse is also true (up to homotopy):
$L(f,g)=0\Longrightarrow$ there are maps $f^{\prime},g^{\prime}$ homotopic to
$f,g$ respectively such that $Coin(f^{\prime},g^{\prime})=\emptyset.$ Now the
problem reads as follows: \textquotedblleft Can we remove coincidences by a
homotopy of $f$ and $g$?\textquotedblright

Let $K=Coin(f,g)$. By $H_{\ast}$ ($H^{\ast}$) we denote the integral singular
(co)homology. For any space $Y$ we define the diagonal map $d:Y\rightarrow
Y\times Y$ by $d(x)=(x,x).$ Let
\[
Y^{\times}=(Y\times Y,Y\times Y\backslash d(Y)).
\]
For codimension $m=0,$ the (cohomology) \textit{coincidence index} $I_{fg}$ of
$(f,g)$ is defined as follows. Since all coincidences lie in $K,$ the map
$(f,g):(N,N\backslash K)\rightarrow M^{\times}$ is well defined. Let $\tau$ be
the generator of $H^{n}(M^{\times})=\mathbf{Z}$ and $O_{N}$ the fundamental
class of $N$ around $K,$ then let
\[
I_{fg}=<(f,g)^{\ast}(\tau),O_{N}>\in\mathbf{Z.}%
\]
The coincidence index satisfies the following natural properties. (1) Homotopy
Invariance: the index is invariant under homotopies of $f,g$; (2) Additivity:
The index over a union of disjoint sets is equal to the sum of the indices
over these sets; (3) Existence of Coincidences: if the index is nonzero then
there is a coincidence; (4) Normalization: the index is equal to the Lefschetz
number; (5) Removability: if the index is zero then the coincidence set can be
removed by a homotopy.

While the coincidence theory for codimension $m=0$ is well developed
\cite[VI.14]{Bredon}, \cite{BS}, \cite{Naka}, \cite[Chapter 7]{Vick}, very
little is known beyond this case. For $m>0,$ the vanishing of the coincidence
index does not always guarantee removability. For codimension $m=$ $1,$ the
secondary obstruction to removability was considered by Fuller
\cite{Fuller-th}, \cite{Fuller} for $M$ simply connected. In the context of
Nielsen Theory the sufficient conditions of the local removability for $m=1$
were studied by Dimovski and Geoghegan \cite{DG}, Dimovski \cite{Dim} for the
projection $f:M\times[0,1]\rightarrow M,$ and by Jezierski \cite{Jez} for
$M,N$ subsets of Euclidean spaces or $M$ parallelizable. Necessary conditions
of the global removability for arbitrary codimension were considered by
Gon\c{c}alves, Jezierski, and Wong \cite[Section 5]{GJW} with $N$ a torus and
$M$ a nilmanifold (see also \cite{GW}).

The main purpose of this note is to provide sufficient conditions of
removability of coincidences for some codimensions higher than $1$. Under a
certain technical condition, the coincidence index defined below is the only
obstruction to removability. This condition holds when (1) $M$ \textit{is a
surface; } (2) \textit{fibers }$f^{-1}(y)$ \textit{of }$f$ \textit{are
acyclic;} or (2) \textit{fibers of }$f$ \textit{are }$m$\textit{-spheres for
}$m=4,5,12$\textit{ and }$n$\textit{ large}. The main theorem partially
complements the results listed above. The proof follows and extends the one of
Brown and Schirmer \cite[Theorem 3.1]{BS} for codimension $0$ (see also Vick
\cite[p. 194]{Vick}).

An area of possible applications is discrete dynamical systems. A
\textit{dynamical system} on a manifold $M$ is determined by a map
$f:M\rightarrow M.$ Then the next position, or state, $f(x)$ depends only on
the current one, $x\in M$. Suppose we have a fiber bundle $F\rightarrow
N^{\underrightarrow{~\ \ \ g\ \ \ ~~}}M$ and a map $f:N\rightarrow M.$ Then
this is a \textit{parametrized\ dynamical system}, where the next position
$f(x,s)$ depends not only on the current one, $x\in M,$ but also the
\textquotedblleft input\textquotedblright, $s\in F.$ Then the Coincidence
Problem asks whether there are a position and an input such that the former
remains unchanged, $f(x,s)=x$. A parametrized dynamical system can also be a
model for a non-autonomous ordinary differential equation: $M$ is the space,
$F$ is the time, and $N$ is the space-time.

\section{Normalization Property.}

For nonzero codimension the homology coincidence index $I_{fg}^{\prime
}=(f,g)_{\ast}(O_{N})$ is replaced with the \textit{homology coincidence
homomorphism }\cite{Brooks1}
\[
I_{fg}^{\prime}=(f,g)_{\ast}:H_{\ast}(N,N\backslash V)\rightarrow H_{\ast
}(M^{\times}),
\]
where $V$ is a neighborhood of $Coin(f,g)$. Let $\pi:M\times M\rightarrow M$
be the projection on the first factor, then $\zeta=(M,\pi,M\times M,d)$ is the
tangent microbundle of $M$ and the Thom isomorphism $\varphi:H_{\ast
}(M^{\times})\rightarrow H_{\ast}(M)$ is given by $\varphi(x)=\pi_{\ast}%
(\tau\frown x),$ where $\tau\in H^{n}(M^{\times}\mathbf{)}$ is the Thom class
of $\zeta.$ The Lefschetz number is replaced with the \textit{Lefschetz
homomorphism} \cite{Sav1} $\Lambda_{fg}:H_{\ast}(N,N\backslash V;\mathbf{Q}%
)\rightarrow H_{\ast}(M;\mathbf{Q})$ of degree $(-n)$ defined as follows.
Suppose $f(N\backslash V)\subset\partial M$. For each $z\in H_{\ast
}(N,N\backslash V),$ let%
\[
f_{!}^{z}=(f^{\ast}D^{-1})\frown z,
\]
where $D:H^{\ast}(M,\partial M;\mathbf{Q})\rightarrow H_{n-\ast}%
(M;\mathbf{Q})$ is the Poincar\'{e}-Lefschetz duality isomorphism
$D(x)=x\frown O_{M}$. Now let
\[
\Lambda_{fg}(z)=\sum_{k}(-1)^{k(k+m)}\sum_{j}x_{j}^{k}\frown g_{\ast}f_{!}%
^{z}(a_{j}^{k}),
\]
where $\{a_{1}^{k},...,a_{m_{k}}^{k}\}$ is a basis for $H_{k}(M)$ and
$\{x_{1}^{k},...,x_{m_{k}}^{k}\}$ the corresponding dual basis for $H^{k}(M).$
Then the Lefschetz-type coincidence theorem \cite[Theorem 6.1]{Sav1} states
that $\varphi I_{fg}^{\prime}=\Lambda_{fg}.$ This is the Normalization
Property, which makes the coincidence homomorphism computable by algebraic means.

Since obstructions to removability of coincidences lie in certain cohomology
groups, we need a cohomological analogue of the theory outlined above. Just as
in the homology case, the cohomology coincidence index can be replaced with
the cohomology coincidence homomorphism.

\begin{definition}
Let $C$ be an isolated subset of $Coin(f,g),$ $W,V$ neighborhoods of $C,$
$C\subset V\subset\overline{V}\subset W\subset N,$ and $W\cap Coin(f,g)=C.$
Then let%
\[
I_{fg}=(f,g)^{\ast}:H^{\ast}(M^{\times})\rightarrow H^{\ast}(W,W\backslash
V).
\]

\end{definition}

However in this paper we consider only the restriction of $I_{fg}$ to
$H^{n}(M^{\times})=\mathbf{Z}$. Therefore the only thing that matters is the
class $I_{fg}(\tau)\in H^{n}(W,W\backslash V),$ where $\tau$ is the generator
of $H^{n}(M^{\times})=\mathbf{Z,}$ which will still be called the
\textit{(cohomology)} \textit{coincidence index}. This index satisfies the
properties of additivity, existence of coincidences and homotopy invariance
proven similarly to Lemmas 7.1, 7.2, 7.4 in \cite[p. 190-191]{Vick} respectively.

We will state the Normalization Property under assumptions similar to the ones
in \cite[Section 2]{Sav}, \cite[Section 5]{Sav1}. Assume that $f(W\backslash
V)\subset\partial M.$

\begin{definition}
For each $z\in H_{n}(W,W\backslash V;\mathbf{Q}),$ define homomorphisms
$\Theta_{q}:H^{q}(M,\partial M;\mathbf{Q})\rightarrow H^{q}(M,\partial
M;\mathbf{Q})$ by%
\[
\Theta_{q}=D^{-1}g_{\ast}(f^{\ast}\frown z).
\]
Then%
\[
L_{z}(f,g)=\sum_{q}(-1)^{q}Tr\Theta_{q}%
\]
is called the (cohomology) \textit{Lefschetz number} with respect to $z$ of
the pair $(f,g).$
\end{definition}

\begin{theorem}
[Normalization]\label{Normalization}Suppose that $f(W\backslash V)\subset
\partial M.$ Then for each $z\in H_{n}(W,W\backslash V;\mathbf{Q}),$%
\[
<I_{fg}(\tau),z>=(-1)^{n}L_{z}(f,g).
\]
Therefore, $L_{z}(f,g)\neq0\Longrightarrow Coin(f,g)\neq\emptyset.$
\end{theorem}

\begin{proof}
The proof repeats the computation in the proof of Theorem 7.12 in \cite[p.
197]{Vick} with Lemmas 7.10 and 7.11 replaced with their generalizations,
Lemmas 3.1 and 3.2 in \cite{Sav}.
\end{proof}

The theorem is true even when $N$ is not a manifold.

\section{Local Removability.}

Let $f:\mathbf{S}^{3}\rightarrow\mathbf{S}^{2}$ be the Hopf map. Then $f$ is
onto, in other words, it has a coincidence with any constant map $c$. However
the coincidence homomorphism $I_{fc}:H^{\ast}((\mathbf{S}^{2})^{\times
})\rightarrow H^{\ast}(\mathbf{S}^{3})$ is zero. Therefore Theorem
\ref{Normalization} fails to detect coincidences. In fact, $f$ has a
coincidence with any map homotopic to $c$ \cite{Brooks}, therefore the
converse of the Lefschetz coincidence theorem for spaces of different
dimensions fails in general. Our main result below is a partial converse.

\begin{theorem}
[Local Removability]\label{Remov}Suppose $f(C)=g(C)=\{u\},u\in M\backslash
\partial M,$ and the following condition is satisfied:%
\[
\text{(A) \ }H^{k+1}(W,W\backslash V;\pi_{k}(\mathbf{S}^{n-1}))=0\text{ for
}k\geq n+1.
\]
Then $I_{fg}(\tau)=0$ implies that $C$ can be removed via a local homotopy of
$f$; specifically, there exists a map $f^{\prime}:N\rightarrow M$ homotopic to
$f$ relative $N\backslash W$ such that
\[
W\cap Coin(f^{\prime},g)=\emptyset.
\]

\end{theorem}

The proof uses the classical obstruction theory. Condition (A)\ guarantees
that only the primary obstruction to the local removability, i.e., the
coincidence index, may be nonzero.

\begin{proof}
We can assume that $U=\mathbf{D}^{n}$ is a neighborhood of $u$ in $M$ such
that $f(W)=U$ and $g(W)\subset U.$ Define $Q:\mathbf{D}^{n}\times
\mathbf{D}^{n}\backslash d(\mathbf{D}^{n})\rightarrow\mathbf{D}^{n}%
\backslash\{0\}$ by $Q(x,y)=1/2(y-x).$ Consider the following commutative
diagram:%
\[%
\begin{array}
[c]{ccccc}%
H^{n-1}(\mathbf{S}^{n-1}) &  &  &  & \\
^{\simeq}\downarrow^{p^{\ast}} &  &  &  & \\
H^{n-1}(\mathbf{D}^{n}\backslash\{0\}) & ^{\underrightarrow{~\ \delta^{\ast
}~\simeq~}} & H^{n}(\mathbf{D}^{n},\mathbf{D}^{n}\backslash\{0\}) &  & \\
\downarrow^{Q^{\ast}} &  & ^{\simeq}\downarrow^{Q^{\ast}} &  & \\
H^{n-1}(\mathbf{D}^{n}\times\mathbf{D}^{n}\backslash d(\mathbf{D}^{n})) &
^{\underrightarrow{~\ \ \ \delta^{\ast}~~}} & H^{n}((\mathbf{D}^{n})^{\times
}) & ^{\underleftarrow{\ ~k^{\ast}~\simeq\ }} & H^{n}(M^{\times})\\
\downarrow^{(f,g)^{\ast}} &  & \downarrow^{(f,g)^{\ast}} & ^{I_{fg}}\swarrow &
\\
H^{n-1}(W\backslash V) & ^{\underrightarrow{~\ \ \ \delta^{\ast}~~}} &
H^{n}(W,W\backslash V). &  &
\end{array}
\]
Here $\delta^{\ast}$ is the connecting homomorphism, $k$ the inclusion, $p$
the radial projection. Let $q=pQ(f,g):W\backslash V\rightarrow\mathbf{S}%
^{n-1}.$ Then $q^{\ast}$ is given in the first column of the diagram.

Now we apply the Extension Theorem, Corollary VII.13.13 in \cite[p.
509]{Bredon}. Suppose $I_{fg}=0.$ Then from the commutativity of the diagram,
$\delta^{\ast}q^{\ast}=0$. Thus the primary obstruction to extending $q$ to
$q^{\prime}:W\rightarrow\mathbf{S}^{n-1},$ $c^{n+1}(q)=\delta^{\ast}q^{\ast},$
vanishes. By condition (A) the other obstructions $c^{k+1}(q),k\geq n,$ also
vanish$.$

Next, $q$ has the form%
\[
q(x)=\dfrac{g(x)-f(x)}{||g(x)-f(x)||}.
\]
Define a map $f^{\prime}:W\rightarrow\mathbf{D}^{n}$ by $f^{\prime
}(x)=g(x)-a(x)q^{\prime}(x),$ where $a:W\rightarrow(0,\infty)$ satisfies the
following: (1) $a$ is small enough so that $f^{\prime}(x)\in\mathbf{D}^{n}$
for all $x\in W,$ (2) $a(x)=||g(x)-f(x)||$ for all $x\in W\backslash V.$ Then
$Coin(f^{\prime},g)=\emptyset$ since $q^{\prime}(x)\neq0.$

To complete the proof observe that $f^{\prime}$ is homotopic to $f|_{W}$
relative $W\backslash V$ because $\mathbf{D}^{n}$ is convex.
\end{proof}

The implications of this result for Nielsen theory will be addressed in a
forthcoming paper.

\section{Further Results.}

Suppose $C=f^{-1}(y),$ where $y\in M\backslash\partial M$ is a regular values
for both $f$ and $f|_{\partial N}.$ Then $C$ is a neat submanifold of $N$ and
it has a tubular neighborhood $T$. Now $T$ can be treated as a disk bundle
$(\mathbf{D}^{m},\mathbf{S}^{m-1})\rightarrow(T,T^{\prime})\rightarrow C.$
Therefore condition (A) takes the form%
\[
\text{(A}^{\prime}\text{) \ }H^{k+1}(T,T^{\prime};\pi_{k}(\mathbf{S}%
^{n-1}))=0\text{ for }k\geq n+1.
\]

In case $C$ is a boundaryless $m$-submanifold of $N,$ we have $H^{n+m}%
(T,T^{\prime};G)=H^{n+m}(T,\partial T;G)=G\oplus...\oplus G.$ Therefore if we
let $k=n+m-1,$ then condition (A$^{\prime}$) implies the following:%
\[
\text{(A}^{\ast}\text{) }\pi_{n+m-1}(\mathbf{S}^{n-1})=0.
\]
This restriction cannot be relaxed, in the following sense. Suppose
\[
\lbrack h]\in\pi_{n+m-1}(\mathbf{S}^{n-1})\backslash\{0\}.
\]
Then $h$ can be extended to a map $f:\mathbf{D}^{n+m}\rightarrow\mathbf{D}%
^{n}\subset M$ by setting $f(0)=0$ and $f(x)=||x||h\left(  \dfrac{x}%
{||x||}\right)  $ for $x\in\mathbf{D}^{n+m}\backslash\{0\}.$ Hence any map
homotopic to $f$ relative $\mathbf{S}^{n+m-1}$ is onto \cite[Theorem VII.5.8,
p. 448]{Bredon}. Therefore coincidences of $f$ and $g,$ where $g$ is constant,
cannot be locally removed.

Below we treat condition (A) as a restriction on an arbitrary fiber
$C=f^{-1}(y),y\in M\backslash\partial M$ of $f.$

\begin{lemma}
\label{Lemma}Suppose for $1\leq p\leq m,$ the submanifold $C$ satisfies
\[%
\begin{array}
[c]{c}%
\text{(1) }H^{p}(C)\otimes\pi_{n+p-1}(\mathbf{S}^{n-1})=0,\text{ and}\\
\text{(2) }Tor(H^{p+1}(C),\pi_{n+p-1}(\mathbf{S}^{n-1}))=0.
\end{array}
\]
Then $C$ satisfies condition (A).
\end{lemma}

\begin{proof}
By the Thom Isomorphism Theorem \cite[Section VI.11]{Bredon}, we have
\[
H^{k+1}(T,T^{\prime};\pi_{k}(\mathbf{S}^{n-1}))=H^{k+1-n}(C;\pi_{k}%
(\mathbf{S}^{n-1})).
\]
By condition (2) and the Universal Coefficient Theorem, Corollary 25.14 in
\cite[p. 263]{Gray}, we have also%
\[
H^{p}(C;\pi_{k}(\mathbf{S}^{n-1}))=H^{p}(C)\otimes\pi_{k}(\mathbf{S}^{n-1}).
\]

\end{proof}

It is known \cite{Toda} that $\pi_{n+m-1}(\mathbf{S}^{n-1})=0,$ for the
following values of $m$ and $n$:

\qquad(1) $m=4$ and $n\geq6;$

\qquad(2) $m=5$ and $n\geq7;$

\qquad(3) $m=12$ and $n=7,8,9,14,15,16,....$

\begin{corollary}
The conclusion of Theorem \ref{Remov} holds when (a) $M$ is a surface; (b)
fibers of $f$ are acyclic; or (c) fibers of $f$ are unions of homology
$m$-spheres for the above values of $m$ and $n$.
\end{corollary}

\begin{proof}
(a) $n=2$ and $\pi_{n+p-1}(\mathbf{S}^{n-1})=0$ for all $p>0.$ (b)
$H^{p}(C)=0$ for all $p>0.$ (c) Either $\pi_{n+p-1}(\mathbf{S}^{n-1})=0$ or
$H^{p}(C)=H^{p+1}(C)=0$ for all $p>0.$ Thus the two conditions of Lemma
\ref{Lemma} are satisfied. Now the conclusion follows from Theorem \ref{Remov}.
\end{proof}

\begin{corollary}
Let $F\rightarrow N^{n+m\underrightarrow{~\ \ \ g\ \ \ ~~}}M^{n}$ be an
$m$-sphere bundle with the above values of $m$ and $n,$ or an $m$-disk bundle.
Then the set $C$ of stationary points of the parametrized dynamical system
$F\rightarrow N^{n+m\underrightarrow{~\ \ \ f,g\ \ \ ~~}}M^{n}$ can be removed
via a local homotopy of $f$ provided $I_{fg}=0.$
\end{corollary}

\section{Coincidence-Producing Maps.}

A boundary preserving map $f:(N,\partial N)\rightarrow(M,\partial M)$ is
called \textit{coincidence-producing} if every map $g:N\longrightarrow M$ has
a coincidence with $f$. Brown and Schirmer \cite[Theorem 7.1]{BS} showed that
if $M$ is acyclic, $\dim N=\dim M=n\geq2,$ then $f$ is coincidence-producing
if and only if $f_{\ast}:H_{n}(N,\partial N)\rightarrow H_{n}(M,\partial M)$
is nonzero. Based on the Normalization and Removability Properties we prove a
generalization of this theorem. We call a map $f:(N,\partial N)\rightarrow
(M,\partial M)$ \textit{weakly coincidence-producing} \cite[Section 5]{Sav} if
every map $g:N\rightarrow M$ with $g_{\ast}=0$ (in reduced homology) has a
coincidence with $f.$ In particular every weakly coincidence-producing map is onto.

A corollary to the Lefschetz type coincidence theorem \cite[Corollary
5.1]{Sav} states that if $f_{\ast}:H_{n}(N,\partial N)\rightarrow
H_{n}(M,\partial M)$ is nonzero then the appropriate Lefschetz homomorphism is
nontrivial and, therefore, $f$ is weakly coincidence-producing. For the
converse we need condition (A)\ as an additional assumption.

\begin{theorem}
Suppose $f$ is boundary preserving and suppose that each fiber $C$ of $f$
satisfies condition (A). Then the following are equivalent:

(1) $f$ is weakly coincidence-producing;

(2) $f_{\ast}:H_{n}(N,\partial N)\rightarrow H_{n}(M,\partial M)$ is nonzero.
\end{theorem}

\begin{proof}
Suppose $f_{\ast}:H_{n}(N,\partial N)\rightarrow H_{n}(M,\partial M)$ is zero.
Choose $g$ to be identically equal to $y\in M\backslash\partial M.$ Then
$C=Coin(f,g)=f^{-1}(y)\subset N\backslash\partial N$. Recall $I_{fg}%
=(f,g)^{\ast}:H^{\ast}(M^{\times})\rightarrow H^{\ast}(N,\partial N).$ Then
for all $z\in H_{n}(N,\partial N),$ we have the following.
\[%
\begin{tabular}
[c]{lll}%
$<I_{fg}^{N}(\tau),z>$ & $=(-1)^{n}L_{z}(f,g),$ by Theorem \ref{Normalization}
& \\
& $=(-1)^{n}Tr\Theta_{n},$ because $g_{\ast}=0$ & \\
& $=(-1)^{n}<f^{\ast}(\overline{O}_{M}),z>,$ where $\overline{O}_{M}$ is the
dual of $O_{M}$ & \\
& $=(-1)^{n}<\overline{O}_{M},f_{\ast}(z)>$ & \\
& $=0.$ &
\end{tabular}
\ \ \ \ \ \ \ \ \ \
\]
Hence $I_{fg}(\tau)=0.$ Therefore by Theorem \ref{Remov} the coincidence set
can be removed. Thus $f$ is not weakly coincidence-producing
\end{proof}

Condition (A) is clearly satisfied for $m=0$. Therefore Brown and Schirmer's
Theorem \cite[Theorem 7.1]{BS} follows. Our theorem also includes the
well-known fact that \textit{a map has degree }$0$\textit{ if and only if it
can be deformed into a map which is not onto}.

Examples of maps satisfying condition (1) of the theorem can be found in
\cite[Section 7]{BS}, see also \cite[Section 6]{Sav1}.

I would like to thank the referee for a number helpful suggestions.

\end{document}